\documentclass[10pt]{article}
\usepackage{}
\usepackage{amssymb}
\usepackage{mathrsfs}
\usepackage{amsmath,amssymb,theorem}
\usepackage{graphicx}
\usepackage{subfigure}
\usepackage{psfrag}
\usepackage{color}
\usepackage{enumerate}
\usepackage{setspace}
\newtheorem{thm}{Theorem}[section]

\newtheorem{lem}[thm]{Lemma}
\theorembodyfont{\rmfamily}
\def\pf{\bigskip\noindent {\bf Proof.}~~}

\def\pf{\bigskip\noindent {\emph{Proof.}}~~}

\def\qed{ \hfill $\square$}

\title{Improved bounds for anti-Ramsey numbers of matchings in outerplanar graphs}

\author{Yifan Pei$^1$, Yongxin Lan$^1$\footnote{Corresponding author}\,\, and Hua He$^1$ \\
$^1$ School of Science\\
Hebei University of Technology, Tianjin, 300401, China\\
Email: imbacheck@163.com; yxlan@hebut.edu.cn; hehua@hehut.edu.cn }
\date{}
\begin{document}\maketitle

\begin{abstract}
Let $\mathcal{O}_n$ be the set of all maximal outerplanar graphs of order $n$. Let $ar(\mathcal{O}_n,F)$ denote the maximum positive integer $k$ such that $T\in \mathcal{O}_n$ has no rainbow subgraph $F$ under a $k$-edge-coloring of $T$. Denote by $M_k$ a matching of size $k$. In this paper, we prove that $ar(\mathcal{O}_n,M_k)\le n+4k-9$ for $n\ge3k-3$, which expressively improves the existing upper bound for $ar(\mathcal{O}_n,M_k)$. We also prove that $ar(\mathcal{O}_n,M_5)=n+4$ for all $n\ge 15$.\\[2mm]
\end{abstract}

\section{Introduction}


 For graphs $H$ and $G$, the {\it anti-Ramsey number of $H$ in $G$}, denoted $ar(G,H)$, is the maximum number of colors in an edge-coloring of $G$ containing no rainbow copy of $H$ (i.e., no two edges of $H$ receive the same colors). The problem of determining anti-Ramsey number of $H$ in $K_n$ dates back to the 1975. 
 Since then, the study of $ar(K_n,H)$ for some special graph $H$ has attracted a lot of attention. Simultaneously, the problem was extended to $ar(G,H)$ when replacing $K_n$ by other graph $G$. Various results on $ar(G,H)$ refer to the survey. 
  Specially, when $H$ is a matching, the value of $ar(G,H)$ has been studied  when $G$ is complete graph, 
  complete split graph  \cite{18JYSC}, complete bipartite graph, 
 regular bipartite graph \cite{17JNWW},  hypergraph 
 and so on.

Let $H$ be a graph and  $\mathcal{G}$ a class of graphs containing at least one graph $G$ such that $H\subseteq G$. The {\it anti-Ramsey number of $H$ in $\mathcal{G}$}, denoted $ar(\mathcal{G},H)$, is the maximum number of colors in an edge-coloring of $G\in \mathcal{G}$ such that $G$ contains no rainbow  $H$ as a subgraph. Hor\v{n}\'ak 
et al.~
began the study of $ar(\mathcal{T}_n,H)$, where $\mathcal{T}_n$ is the class of all plane triangulations of order $n$. They considered the anti-Ramsey numbers of cycles in $\mathcal{T}_n$, which were improved in \cite{19LSS}.
For large $n$, the authors determined the bounds for the anti-Ramsey numbers of matchings in $\mathcal{T}_n$.
Let $\mathcal{O}_n$ denote the class of all maximal outerplanar graphs of order $n$.  Jin and Ye \cite{18JY}  investigated the anti-Ramsey numbers of matchings in $\mathcal{O}_n$. They determined the exact values of $ar(\mathcal{O}_n,M_k)$ for all $2\le k\le 4$ and gave the upper and lower bounds for general matchings. We summarize their results as follows.

\begin{thm}[\cite{18JY}]\label{known}
Let $n$ and $k$ be positive integers. Then
\begin{description}
\item[(1)] $ar(\mathcal {O}_n, M_2)=1$ for $n\ge5$, and $ar(\mathcal {O}_4, M_2)=3$;

\item[(2)] $ar(\mathcal {O}_n, M_3)=n$ for $n\ge7$, and $ar(\mathcal {O}_6, M_3)=7$;

\item[(3)] $ar(\mathcal {O}_n, M_4)=n+2$ for $n\ge9$, and $ar(\mathcal {O}_8, M_4)=11$;

\item[(4)] $n+2k-6 \le ar(\mathcal {O}_n, M_k)\le n+14k-25$ for all $k \ge 5$ and $n\ge 2k$.

%
\end{description}
\end{thm}

In this paper, we first improve the existing upper bound for $ar(\mathcal {O}_{n}, M_k)$.
\begin{thm}\label{upper}
For all $n\ge 3k-3$ and $k\ge2$, $ar(\mathcal {O}_{n}, M_k)\le n+4k-9$.
\end{thm}

Theorem~\ref{upper} expressively improves the upper bound in Theorem~\ref{known}(4). We believe that the lower bound obtained  in Theorem~\ref{known}(4) is the exact value of  $ar(\mathcal {O}_{n}, M_k)$ for large $n$. This is indeed the case for $M_5$.

\begin{thm}\label{main}
For all $n\ge 15$, $ar(\mathcal {O}_{n}, M_5)=n+4$.
\end{thm}

The following famous Tutte-Berge formula will be very helpful in the proof of Theorem~\ref{main}. Let $\beta(G)$ be the matching number of $G$.

\begin{thm}
\label{TBF}
 Let $G$ be a graph of order $n$. Then there exists  $T\subseteq V(G)$ with $|T|\le m(G)$ such that
$$\beta(G)=\frac{1}{2}(n-o(G-T)+|T|),$$
where $o(G-T)$ is the number of odd components in $G-T$. Moreover, each odd component of $G-T$ is factor-critical and each even component of $G-T$ has a perfect matching.
\end{thm}

In order to prove our results, we first show that the following significant result.

\begin{lem}\label{bipartite}
Let $G$ be a bipartite outerplanar graph of order $n$ with two parts $X$ and $Y$ and $|Y|\ge|X|\ge1$. Then $e(G)\le n+|X|-2$.\vspace{-4mm}
\end{lem}
\pf Let $G$ be given as in the statement. We shall prove that $e(G)\le n+|X|-2$ by induction on $n$.
It can be easily checked that the result holds when $n=2$. So we assume that $n\ge3$. Since $G$ is an outerplanar graph, we see that $\delta(G)\le 2$.  Let $v\in V(G)$ with $d_G(v)=\delta(G)$. Assume that $d_G(v)\le 1$. By induction hypothesis, $e(G-v)\le n+|X|-3$ and so $e(G)=e(G-v)+d_G(v)\le n+|X|-2$, as desired. So we next assume that $\delta(G)=2$. We first consider $v\in X$.
By induction hypothesis, $e(G-v)\le n-1+|X-v|-2$ and so $e(G)=e(G-v)+d_G(v)\le n+|X|-2$. So we  further assume that $v\in Y$.
Let $N_G(v)=\{v_1,v_2\}$. By induction hypothesis, $e(G-v)\le n+|X|-3$.
 Suppose $e(G-v)=n+|X|-3$. Then $G-v$ is a maximal bipartite outerplanar graph. Hence, the size of each face of $G-v$ except for the outer face is four. Then $G-v$ contains a $K_{2,2}$-minor with one part $\{v_1, S\}$, where $G[S]$ is a connected subgraph containing $v_2$. But then $G$ has a $K_{2,3}$-minor, which contradicts to the outerplanarity of $G$. Thus, $e(G-v)\le n+|X|-4$ and so $e(G)=e(G-v)+d_G(v)\le n+|X|-2$, as desired. \qed
\medskip
\section{Proof of Theorem~\ref{upper}}

We process it by induction on $k$. By Theorem~\ref{known}(1-3), the results hold when $k\le4$. So we next assume $k\ge5$. Suppose $ar(\mathcal {O}_{n}, M_k)\ge n+4k-8$ for some $k$ and $n\ge3k-3$. Then there must have $G\in \mathcal {O}_{n}$ such that there is no rainbow subgraph $M_k$ in some edge-colored graph $G$ with at least $n+4k-8$ colors. By induction hypothesis, $G$ contains a rainbow copy of $M_{k-1}$. Let $H$ be a rainbow spanning subgraph of $G$ containing $M_{k-1}$ as a subgraph. Then $e(H)\ge n+4k-8$. Let $M$ be a matching of size $k-1$ with vertex set $\{v_1,v_2,\ldots, v_{2(k-1)}\}$ and edge set $\{v_{2i-1}v_{2i}: 1\le i\le k-1\}$. Let $I=V(H)-V(M)$.  Let $M'$ be a matching of size $t$ consisting of edges $e\in M$ such that both ends of $e$ have neighbour in $I$. Then for any edge $f\in M-M'$, at least one end of $f$ has no neighbour in $I$. Without loss of generality, we may assume $M'=\{v_1v_2,v_3v_4,\ldots,v_{2t-1}v_{2t}\}$ and $d_I(v_{2i-1})=0$ for any $t+1\le i\le k-1$. Since $G$ has no rainbow subgraph $M_k$, it follows that $I$ is an independent set in $H$ and $d_{I}(u)=1$ for any $u\in V(M')$. Then $e_{H}(V(M'),I)=2t$. Let $S=\{v_{2(t+1)},v_{2(t+2)},\ldots,v_{2(k-1)}\}$. By Lemma~\ref{bipartite}, $e_H(S,I)\le 2|S|+|I|-2=2(k-1-t)+(n-2(k-1))-2=n-2t-2$ because $|I|=n-2(k-1)\ge k-1\ge|S|$. Since $G$ is an outerplanar graph, we see $e(H[V(M)])\le 2|V(M)|-3=4(k-1)-3=4k-7$. Thus, $e(H)=e(H[V(M)])+e_{H}(V(M'),I)+e_H(S,I)\le (4k-7)+2t+(n-2t-2)=n+4k-9$, a contradiction. Hence, $ar(\mathcal {O}_{n}, M_k)\le n+4k-9$ for any $k\ge2$ and $n\ge3k-3$. This completes the proof of Theorem~\ref{upper}.\qed

\section{Proof of Theorem \ref{main}.}
%
By Theorem~\ref{known}(4), $ar(\mathcal{O}_n,M_5)\ge n+4$. We shall prove that   $ar(\mathcal{O}_n,M_5)\le n+4$.
Suppose not.  Then there is a $G\in \mathcal{O}_n$  satisfying $G$ contains no rainbow $M_5$ under some surjection $\alpha: E(G)\rightarrow [k]$, where $k\ge n+5$.
 We use $H$ to denote a rainbow spanning subgraph of $G$ such that $e(H)=k$. By Theorem~\ref{known}(3), $H$ contains $M_4$ as a subgraph. Clearly, $\beta(H)=4$ because $G$ contains no rainbow $M_5$. By Theorem~\ref{TBF}, there is a $T\subseteq V(H)$ with $t:=|T|\le \beta(H)=4$ satisfying $p:=o(H-T)=n+t-8$.
Let $A_1,\ldots, A_p$ denote all odd components of $H-T$ with $|A_i|=a_i$ for each $i\in[p]$.  We may assume that $a_1\le \cdots \le a_p$. Let $r:=\max\{i: a_i=1\}$.
Let $V(A_i)=\{u_i\}$ for all $i\in[r]$  with $d_H(u_1)\ge \cdots \ge d_H(u_r)$. Moreover, let $T:=\{w_1,\ldots,w_t\}$ when $t\ge1$, and $U:=\{u_1,\ldots,u_r\}$. Let $R$ denote the set of vertices of all even components of $H-T$.
We claim $t\ge1$. Suppose $t=0$. Then $p=n-8$. Hence, $e(H)=\sum_{i=r+1}^{p-1}e(A_i)+e(H[V(A_{p})\cup R])\le \sum_{i=r+1}^{p-1}(2|A_{i}|-3)+2|V(A_{p})\cup R|-3= 2(n-r)-3(p-r)=2(n-r)-3(n-8-r)=24-(n-r)\le 16<n+5$ because $r\le p=n-8$ and $n\ge 15$, a contradiction. Thus, $t\ge1$. We shall prove several useful claims. \\

\noindent {\bf Claim 1.} {\it If $H-\{x,y\}$ has two edge-disjoint $M_4$, then $xy\notin E(G)$.}\medskip

\noindent {\bf Proof.} Suppose not. Let $M_1$ and $M_2$ be two edge-disjoint $M_4$ in $H-\{x,y\}$. If $\alpha(xy)=\alpha(e)$ for some $e\in M_1$, then $M_2\cup \{xy\}$ is a rainbow $M_5$ in $G$, a contradiction. Thus, $\alpha(xy)\neq \alpha(e)$ for all $e\in M_1$. But then $M_1\cup \{xy\}$ is a rainbow $M_5$ in $G$, a contradiction.\qed\\

\noindent {\bf Claim 2.} {\it Let $xy\in E(G)$ and $uv\in E(H)$. If $H-\{x,y,u,v\}$ has two edge-disjoint $M_3$, then $\alpha(xy)=\alpha(uv)$.}\medskip

\noindent {\bf Proof.} Suppose not. Let $M_1$ and $M_2$ be two edge-disjoint $M_3$ in $H-\{x,y,u,v\}$. If $\alpha(xy)=\alpha(e)$ for some $e\in M_1$, then $M_2\cup\{xy,uv\}$ is a rainbow $M_5$ in $G$, a contradiction. Thus, $\alpha(xy)\neq \alpha(e)$ for all $e\in M_1$. But then $M_1\cup\{xy,uv\}$ is a rainbow $M_5$ in $G$, a contradiction.\qed\\


\noindent {\bf Claim 3.} $R=\emptyset$.\medskip

\noindent{\bf Proof.} Suppose not. That is, $|R|\ge2$. Note that $|R|\le n-p-t=8-2t$. Then $t\le3$. Assume $t=3$. Then $|R|=2$ and $a_i=1$ for each $i\ge1$. By Lemma~\ref{bipartite}, $e_H(T,U\cup R)\le n+1$. Hence, $n+5\le e(H)=e_H(T,U\cup R)+e(H[T])+e(H[R])\le (n+1)+3+1=n+5$. It follows that $H[T]=K_3$ and $e_H(T,U\cup R)=n+1$, which implies that there exist four vertices $x_1,x_2,x_3,x_4\in U\cup R$ such that $d_T(x_i)=2$. But then $G$ has a $K_{2,3}$-minor, a contradiction. Assume next $t=2$. Then either $|R|=2$ or $|R|=4$. If $|R|=2$, then $a_p=3$ and $a_i=1$ for each $i\in[p-1]$. By Lemma~\ref{bipartite}, $e_H(T,U\cup R)\le n-3$. Hence, $n+5\le e(H)=e(H[T\cup V(A_p)])+e_H(T,U\cup R)+e(H[R])\le (2\cdot5-3)+(n-3)+1=n+5$, which implies that $H[T\cup V(A_p)]$ is a maximal outerplanar graph and $d_T(x)=d_T(y)=2$ for some $x,y\in U\cup R$. But then $G$ has a $K_{2,3}$-minor (with one part $T$ and the other part $\{x,y,V(A_p)\}$), a contradiction. Thus, $|R|=4$ and so $a_i=1$ for each $i\in[p]$. By Lemma~\ref{bipartite}, $e_H(T,U)\le n-4$. Then $n+5\le e(H)=e_H(T,U)+e(H[T\cup R])\le (n-4)+(2\cdot6-3)=n+5$, which implies that $d_H(u_1)=d_H(u_2)=2$ and $H[T\cup R]$ is a maximal outerplanar graph. But then $G$ contains a $K_{2,3}$-minor (with one part $T$ and the other part $\{u_1,u_2,R\}$), a contradiction.
Thus, $t=1$ and so $p=n-7$. However, $e(H)=\sum_{i=1}^p e(H[T\cup A_i])+e(H[T\cup R])=\sum_{i=1}^p (2(|A_i|+t)-3)+2(|R|+t)-3=2(\sum_{i=1}^p |A_i|+|R|)-(p+1)=2(n-1)-(p+1)=n+4$, a contradiction.\qed\\


\noindent{\bf Claim 4.} $a_p\le3$.\medskip

\noindent{\bf Proof.} Suppose $a_p\ge 5$. By Claim 3, $t\le 2$, else $p=n+t-8\ge n-5$ and so $a_p\le3$. Assume $t=2$. Then  $a_p=5$ and $a_i=1$ for each $i\in [p-1]$. We claim $e(H[T\cup V(A_p)])\le 2\cdot 7-4$. Suppose $e(H[T\cup V(A_p)])= 2\cdot 7-3$. That is $H[T\cup V(A_p)]$ is a maximal outerplanar graph, which yields $d_H(u_i)=1$ for each $i\ge2$. Because $e_H(T,U)=e(H)-e(H[T\cup V(A_p)])\ge n-6$, we see $d_H(u_1)=2$. Note that $H[T\cup V(A_p)\cup\{u_1\}]$ is a maximal outerplanar graph and contains two edge-disjoint $M_4$.  By Claim 1, $E(G[U-u_1])=\emptyset$. Thus, $e_G(u_i,V(A_p)\cup \{u_1\})\ge1$ for each $i\ge2$ because $\delta(G)\ge 2$. Note that $|U|\ge 8$ as $n\ge15$. Then either $d_{U-u_1}(u_1)\ge3$ or $d_{U-u_1}(V(A_p))\ge3$. But then $G$ has a $K_{2,3}$-minor (with one part $\{u_i,u_j,u_k\}\subseteq U-u_1$ and the other part $\{T,u_1\}$ or $\{T,V(A_p)\}$), a contradiction.
Thus, $e(H[T\cup V(A_p)])\le 2\cdot 7-4$. By Lemma~\ref{bipartite}, $e_H(T,U)\le n-5$. Then $n+5\le e(H)=e_H(T,U)+e(H[T\cup V(A_p)])\le (n-5)+(2\cdot7-4)=n+5$, which implies that $e(H[T\cup V(A_p)])=10$, $d_H(u_1)=d_H(u_2)=2$ and $d_H(u_i)=1$ for each $i\ge3$. Thus, $e_H(w_i,V(A_p))=0$ for some $i\in[2]$. Without loss of generality, we may assume $e_H(w_1,V(A_p))=0$. Then $H[V(A_p)\cup \{w_2\}]$ is a maximal outerplanar graph because $e(H[T\cup V(A_p)])=10$. Note that $H-\{u_i,u_j\}$ for $i,j\ge3$ contains two edge-disjoint $M_4$. By Claim 1, $E(G[U-\{u_1,u_2\}])=\emptyset$. Thus, $e_G(u_i,V(A_p)\cup \{u_1,u_2\})\ge1$ for each $i\ge3$ because $\delta(G)\ge 2$. We claim $d_{U-\{u_1,u_2\}}(\{u_1,u_2\})\le 2$. Suppose not. Then $u_1u_i,u_2u_j\in E(G)$ for some $i,j\ge3$. Note that $H[V(A_p)\cup \{w_2\}]$ contains two edge-disjoint $M_3$, say $M$ and $M'$. By Claim 2, $\alpha(u_1u_i)=\alpha(u_2w_1)$ and $\alpha(u_2u_j)=\alpha(u_1w_1)$. But then $G$ has a rainbow $M_5=\{u_1u_i,u_2u_j\}\cup M$, a contradiction. Thus, $d_{U-\{u_1,u_2\}}(\{u_1,u_2\})\le 2$. Note that $|U|\ge 8$ as $n\ge 15$. This implies $d_{U-\{u_1,u_2\}}(V(A_p))\ge 3$. But then $G$ has a $K_{2,3}$-minor (with one part $\{x,y,z\}\subseteq U-\{u_1,u_2\}$ and the other part $\{T\cup\{u_1\}, V(A_p)\}$), a contradiction.

 Thus, $t=1$. Then either $a_p=5$ or $a_p=7$. Suppose $a_p=7$. Then $a_i=1$ for each $i\in[p-1]$. We see $n+5\le e(H)=e_H(T,U)+e(H[T\cup V(A_{p})])\le (n-8)+(2\cdot8-3)=n+5$, which implies $H[T\cup V(A_p)]$ is a maximal outerplanar graph and $d_H(x)=1$ for each $x\in U$. Note that $|U|\ge7$ as $n\ge15$ and $H[T\cup V(A_p)]$ contains two edge-disjoint $M_4$. By Claim 1, $E(G[U])=\emptyset$. Since $\delta(G)\ge 2$, $e_{G}(x,V(A_p))\ge1$ for each $x\in U$. But then $G$ contains a $K_{2,3}$-minor (with one part $\{u_1,u_2,u_3\}$ and the other part $\{w_1,V(A_p)\}$), a contradiction. Thus, $a_p=5$. Then $a_{p-1}=3$ and $a_i=1$ for each $i\in[p-2]$. Hence, $n+5\le e(H)=e_H(T,U)+e(H[T\cup V(A_{p-1})])+e(H[T\cup V(A_{p})])\le (n-9)+(2\cdot4-3)+(2\cdot6-3)=n+5$, which implies  $H[T\cup V(A_{i})]$ for any $i\in\{p-1,p\}$ is a  maximal outerplanar graph and $d_H(x)=1$ for each $x\in U$. Note that $H[T\cup V(A_{p-1})\cup V(A_p)]$ contains two edge-disjoint $M_4$. By Claim 1, $E(G[U])=\emptyset$. Since $\delta(G)\ge 2$, $e_{G}(x,V(A_{p-1})\cup V(A_p))\ge1$ for each $x\in U$. Note that $|U|\ge 6$ as $n\ge 15$. Then $d_{U}(V(A_i))\ge3$ for some $i\in\{p-1,p\}$. But then $G$ has a $K_{2,3}$-minor (with one part $\{u_j,u_k,u_{\ell}\}\subseteq U$ and the other part $\{w_1,V(A_{i})\}$), a contradiction.\qed\\

Now we will continue to complete our proof. We first consider $t=1$. Then $p=n-7$. By Claims 3 and 4,
$a_p=a_{p-1}=a_{p-2}=3$ and $a_i=1$ for each $i\in[p-3]$. By Theorem~\ref{TBF}, $A_p=A_{p-1}=A_{p-2}=K_3$. Hence, $n+5\le e(H)=e_H(T,U)+\sum_{i=p-2}^pe(H[T\cup V(A_{i})])\le (n-10)+\sum_{p-2}^p(2\cdot4-3)=n+5$, which implies $H[T\cup V(A_{i})]$ for each $i\ge p-2$ is a  maximal outerplanar graph and $d_H(x)=1$ for each $x\in U$. Note that $|U|\ge5$ as $n\ge15$ and $H[T\cup V(A_{p-2})\cup V(A_{p-1})\cup V(A_p)]$ contains two edge-disjoint $M_4$. By Claim 1, $E(G[U])=\emptyset$. Since $\delta(G)\ge 2$, $e_{G}(x,V(A_{p-2})\cup V(A_{p-1})\cup V(A_p))\ge1$ for each $x\in U$. 
Then there exist two vertices in $U$, say $u_1,u_2$, such that $e_G(u_1,V(A_i))\ge1$ and $e_G(u_2,V(A_i))\ge1$ for some $i\ge p-2$. Without loss of generality, we may assume $e_G(u_i,V(A_p))\ge1$ for each $i\in[2]$. Let $V(A_{p})=\{v_1,v_2,v_3\}$. Assume $u_1v_1,u_2v_2\in E(G)$. Note that $H[T\cup V(A_{p-2})\cup V(A_{p-1})]$ contains two edge-disjoint $M_3$, say $M$ and $M'$. By Claim 2, $\alpha(u_1v_1)=\alpha(v_2v_3)$ and $\alpha(u_2v_2)=\alpha(v_1v_3)$. But then $G$ has a rainbow $M_5=\{u_1v_1,u_2v_2\}\cup M$, a contradiction.

\medskip
Assume next $t=2$. Then $p=n-6$.
By Claims 3 and 4,
$a_p=a_{p-1}=3$ and $a_{i}=1$ for each $i\in[p-2]$. By Theorem~\ref{TBF}, $A_p=A_{p-1}=K_3$.
Let $V(A_{p-1})=\{v_1,v_2,v_3\}$ and $V(A_{p})=\{v_4,v_5,v_6\}$. Note that $|U|\ge 7$ as $n\ge 15$.
We first claim $d_H(u_1)=2$. Suppose not. Then $n+5\le e(H)=e_H(T,U)+e(H[T\cup V(A_{p-1})\cup V(A_p)])\le (n-8)+(2\cdot 8-3)=n+5$, which implies that $H[T\cup V(A_{p-1})\cup V(A_p)]$ is a maximal outerplanar graph and so contains two edge-disjoint $M_4$. By Claim 1, $E(G[U])=\emptyset$. Thus, $e_G(u_i,V(A_{p-1})\cup V(A_p))\ge 1$ for each $i\ge1$ because $\delta(G)\ge 2$.  Hence, either $d_{U}(V(A_{p-1}))\ge 3$ or $d_{U}(V(A_{p}))\ge 3$ because $|U|\ge7$. But then $G$ has a $K_{2,3}$-minor (with one part $\{x,y,z\}\subseteq U$ and the other part $\{T, V(A_{p-1})\}$ or $\{T, V(A_{p})\}$), a contradiction.
 Thus, $d_H(u_1)=2$. Then $e(H[T\cup V(A_{p-1})\cup V(A_p)])\le 2\cdot 8-4$. We next claim $d_H(u_2)=2$. Suppose not. Then $n+5\le e(H)=e_H(T,U)+e(H[T\cup V(A_{p-1})\cup V(A_p)])\le (n-7)+12=n+5$, which implies that $e(H[T\cup V(A_{p-1})\cup V(A_p)])=12$. Without loss of generality, assume $w_1v_1,w_2v_2,w_2v_4,w_2v_5\in E(H)$. Note that $H-\{u_i,u_j\}$ contains two edge-disjoint $M_4$ for any $i,j\ge2$. By Claim 1, $E(G[U-u_1])=\emptyset$. Thus, $e_G(u_i,V(A_{p-1})\cup V(A_p)\cup\{u_1\})\ge 1$ for each $i\ge 2$ because $\delta(G)\ge 2$.
We claim $d_{U-u_1}(V(A_p))\le1$. Suppose not. Assume $u_iv_4,u_jv_5\in E(G)$. Note that $H[T\cup V_{p-1}\cup\{u_1\}]$ contains two edge-disjoint $M_3$, say $M$ and $M'$. By Claim 2, $\alpha(u_iv_4)=\alpha(v_5v_6)$ and $\alpha(u_jv_5)=\alpha(v_4v_6)$. But then $M\cup\{u_iv_4,u_jv_5\}$ is a rainbow $M_5$, a contradiction, as desired.  Then either $d_{U-u_1}(V(A_{p-1}))\ge3$ or $d_{U-u_1}(u_1))\ge3$ because $|U|\ge7$. But then $G$ has a $K_{2,3}$-minor (with one part $\{x,y,z\}\subseteq U-u_1$ and the other part $\{T, V(A_{p-1})\}$ or $\{T, \{u_1\}\}$), a contradiction.
 Thus, $d_H(u_2)=2$. That is $d_H(u_i)\le1$ for any $i\ge3$ which means $e_H(T,U)\le n-6$. Hence,  $e_H(T,V(A_i))\le2$ for each $i\in\{p-1,p\}$ and so $e(H[T\cup V(A_{p-1})\cup V(A_p)])\le 11$. Thus, $n+5\le e(H)=e_H(T,U)+e(H[T\cup V(A_{p-1})\cup V(A_p)])\le (n-6)+11=n+5$, which implies that $e_H(T,V(A_i))=2$ for each $i\in\{p-1,p\}$ and $w_1w_2\in E(H)$. We claim $N_T(V(A_{p-1}))=N_{T}(V(A_{p}))$. Suppose not. Then $H-\{u_i,u_j\}$ contains two edge-disjoint $M_4$ for $i,j\ge1$. By Claim 1, $E(G[U])=\emptyset$. Thus, $e_G(u_i,V(A_{p-1})\cup V(A_p))\ge1$ for each $i\ge3$ because $\delta(G)\ge2$.
 Then $d_{U-\{u_1,u_2\}}(V(A_i))\ge3$ for some $i\in\{p-1,p\}$ because $|U|\ge7$. Similarly, $G$ has a $K_{2,3}$-minor, a contradiction. Thus, assume $N_T(V(A_{p-1}))=N_{T}(V(A_{p}))=\{w_2\}$. Note that $H-\{u_i,u_j\}$ contains two edge-disjoint $M_4$ for any $i,j\ge3$.
 By Claim 1, $E(G[U-\{u_1,u_2\}])=\emptyset$. Thus, $e_G(u_i,V(A_{p-1})\cup V(A_p)\cup\{u_1,u_2\})\ge1$ for each $i\ge3$.
 Since $H-V(A_i)$ for any $i\in\{p-1,p\}$ and $H-\{u_1,u_2,w_1\}$ contain two edge-disjoint $M_3$, we see that $e_G(V(A_i),U-\{u_1,u_2\})\le1$ and $e_G(\{u_1,u_2\},U-\{u_1,u_2\})\le2$. Hence, $|U-\{u_1,u_2\}|\le 4$ and so $|U|\le 6$, which contradicts to $n\ge15$.\medskip

Assume then $t=3$. Then $p=n-5$.
By Claim 3, $a_p=3$ and $a_i=1$ for each $i\in[p-1]$. By Theorem~\ref{TBF}, $A_p=K_3$. Note that $|U|\ge9$ as $n\ge15$. By Lemma~\ref{bipartite}, $e_H(T,U)\le n-2$. We see $e_H(T,U)\ge n-4$, else $e(H)=e_H(T,U)+e(H[T\cup V(A_p)])\le n+4$. Thus, $e_H(T,U)\in\{n-2,n-3,n-4\}$.

Suppose $e_H(T,U)=n-2$. Then either $d_H(u_1)=3$ and $d_H(u_2)=d_H(u_3)=2$ or $d_H(u_i)=2$ for all $i\in[4]$. It implies $e(H[T])\le 2$ and $e_H(T,V(A_p))\le 2$. So $n+5\le e(H)=e_H(T,U)+e(H[T])+e_H(T,V(A_p))+e(A_p)\le n+5$, which means $e(H[T])=2$ and $e_H(T,V(A_p))=2$. For some $\ell\in\{4,5\}$, $H-\{u_i,u_j\}$ contains two edge-disjoint $M_4$ for any $i,j\ge \ell$. By Claim 1, $E(G[U-\{u_1,\ldots,u_{\ell-1}\}])=\emptyset$. Hence, $e_G(u_i,V(A_p)\cup \{u_1,\ldots,u_{\ell-1}\})\ge 1$ for each $i\ge \ell$ because $\delta(G)\ge2$. It is easy to see $d_{U-\{u_1,\ldots,u_{\ell-1}\}}(\{u_1,\ldots,u_{\ell-1}\})\le 2$.  Thus, $d_{U-\{u_1,\ldots,u_{\ell-1}\}}(V(A_p))\ge 3$ because $|U|\ge9$. Obviously, $G$ has a $K_{2,3}$-minor, a contradiction.

 Suppose $e_H(T,U)=n-3$. Then either $d_H(u_1)=3$ and $d_H(u_2)=2$ or $d_H(u_i)=2$ for all $i\in[3]$. It implies $e(H[T])+e_H(T,V(A_p))\le 5$  and $e_H(T,V(A_p))\le3$. Hence, $n+5\le e(H)=e_H(T,U)+e(H[T])+e_H(T,V(A_p))+e(A_p)\le n+5$, which means $e(H[T])+e_H(T,V(A_p))=5$. We claim $d_H(u_i)=2$ for all $i\in[3]$. Suppose not. Note that $H-\{u_i,u_j\}$ contains two edge-disjoint $M_4$ for any $i,j\ge3$. By Claim 1, $E(G[U-\{u_1,u_2\}])=\emptyset$. Hence, $e_G(u_i,V(A_p)\cup \{u_1,u_2\})\ge 1$ for each $i\ge 3$ because $\delta(G)\ge2$. Then either $d_{U-\{u_1,u_2\}}(u_i)\ge 3$ for some $i\in[2]$ or $d_{U-\{u_1,u_2\}}(V(A_p))\ge 3$ because $|U|\ge9$, which means $G$ has a $K_{2,3}$-minor, a contradiction.
Thus, $d_H(u_i)=2$ for all $i\in[3]$. Note that $H-\{u_i,u_j\}$ contains two edge-disjoint $M_4$ for any $i,j\ge4$. By Claim 1, $E(G[U-\{u_1,u_2,u_3\}])=\emptyset$. Hence, $e_G(u_i,V(A_p)\cup \{u_1,u_2,u_3\})\ge 1$ for each $i\ge 4$ because $\delta(G)\ge2$. It is easy to see $d_{U-\{u_1,u_2,u_3\}}(\{u_1,u_2,u_3\})\le 3$. Thus, $d_{U-\{u_1,u_2,u_3\}}(V(A_p))\ge 3$ because $|U|\ge9$. Obviously, $G$ has a $K_{2,3}$-minor, a contradiction.

Suppose $e_H(T,U)=n-4$.  Then either $d_H(u_1)=3$ or $d_H(u_i)=2$ for all $i\in[2]$. If $d_H(u_1)=3$, then $e(H[T])\le 2$ and $e_H(T,V(A_p))\le 3$. Hence, $e(H)=e_H(T,U)+e(H[T])+e_H(T,V(A_p))+e(A_p)\le n+4$, a contradiction. Thus, $d_H(u_i)=2$ for all $i\in[2]$. Then $N_H(u_1)\neq N_H(u_2)$, else $e(H)=e_H(T,U)+e(H[T])+e_H(T,V(A_p))+e(A_p)\le (n-4)+2+3+3=n+4$. Let $N_H(u_1)=\{w_1,w_2\}$ and $N_H(u_2)=\{w_2,w_3\}$. Since $n+5\le e(H)=e_H(T,U)+e(H[T\cup V(A_p)])\le (n-4)+9=n+5$, we see $e(H[T\cup V(A_p)])$ is a maximal outerplanar graph, which means $w_1x,w_3y\in E(H)$ for some $x\neq y\in V(A_p)$ and $e(H[T])\ge2$. Note that $H-\{u_i,u_j\}$ contains two edge-disjoint $M_4$ for any $i,j\ge3$. By Claim 1, $E(G[U-\{u_1,u_2\}])=\emptyset$. Hence, $e_G(u_i,V(A_p)\cup \{u_1,u_2\})\ge 1$ for each $i\ge 3$ because $\delta(G)\ge2$. Thus, either $d_{U-\{u_1,u_2\}}(u_i)\ge 3$ for some $i\in[2]$ or $d_{U-\{u_1,u_2\}}(V(A_p))\ge 3$ because $|U|\ge9$, which means $G$ has a $K_{2,3}$-minor, a contradiction.\medskip

Finally, assume $t=4$. Then $p=n-4$, which means $|R|=0$ and $a_i=1$ for each $i\in[p]$. That is $|U|=p=n-4\ge11$ because $n\ge15$. By Lemma~\ref{bipartite}, $e_H(T,U)\le n+2$. We see $e_H(T,U)\ge n$, else $e(H)=e_H(T,U)+e(H[T])\le n+4$. Thus, $e_H(T,U)\in\{n+2,n+1,n\}$.

Assume $e_H(T,U)=n+2$. Then $e(H[T])=e(H)-e_H(T,U)\ge3$. We next prove that \medskip

\noindent{($\ast$) $d_H(u_1)+d_H(u_2)\le5$.}\medskip

To see why ($\ast$) holds, suppose $d_H(u_1)+d_H(u_2)\ge6$. Then either $d_H(u_1)=4$ and $d_H(u_2)=2$ or $d_H(u_1)=d_H(u_2)=3$. Hence, $d_H(u_3)=d_H(u_4)=2$ and $d_H(u_i)=1$ for each $i\ge5$. Note that $H-\{u_i,u_j\}$ contains two edge-disjoint $M_4$ for any $i,j\ge5$. By Claim 1, $E(G[U-\{u_1,\ldots,u_4\}])=\emptyset$. Hence, $e_G(u_i,\{u_1,\ldots,u_4\})\ge 1$ for each $i\ge 5$ because $\delta(G)\ge2$. Note that $|U|\ge11$. Then there exist two adjacent vertices in $T$, say $w_1,w_{2}$, such that $d_{U-\{u_1,\ldots,u_4\}}(\{w_1,w_2\})\ge3$ and $N_H(u_i)=\{w_1,w_{2}\}$ for some $i\in\{3,4\}$.
Assume $u_5,u_6\in N_{U-\{u_1,\ldots,u_4\}}(w_1)$  and $u_7\in N_{U-\{u_1,\ldots,u_4\}}(w_2)$. Then either $u_{i}u_5\in E(G)$ or $u_{i}u_6\in E(G)$.
Note that both $H-\{u_{i},u_5\}$ and  $H-\{u_{i},u_6\}$ contain two edge-disjoint $M_4$. But then $u_{i}u_5\notin E(G)$ and $u_{i}u_6\notin E(G)$ by Claim 1, a contradiction. Thus, $d_H(u_1)+d_H(u_2)\le5$. This proves $(\ast)$. \medskip

Note that $d_H(u_2)\ge2$. By $(\ast)$, either $d_H(u_1)=3$ and $d_H(u_2)=\cdots=d_H(u_5)=2$
 and $d_H(u_i)=1$ for each $i\ge 6$, or $d_H(u_1)=\cdots=d_H(u_6)=2$ and $d_H(u_i)=1$ for each $i\ge7$. Let $N_H(u_2)=\{w_1,w_2\}\subseteq N_H(u_1)$, $N_H(u_3)=N_H(u_4)=\{w_3,w_4\}$ and $N_H(u_6)\subseteq N_H(u_5)=\{w_2,w_3\}$.
 Note that $H-\{u_i,u_j\}$ contains two edge-disjoint $M_4$ for any $i,j\ge5$. By Claim 1, $E(G[U-\{u_1,\ldots,u_4\}])=\emptyset$. Hence, $e_G(u_i,\{u_1,\ldots,u_4\})\ge 1$ for each $i\ge 7$. By Claims 1 and 2, we can easily see that $d_{U-\{u_1,\ldots,u_6\}}(u_i)\le1$  for any $i\in[4]$. Hence, $|U-\{u_1,\ldots,u_6\}|\le4$, which contradicts to $|U|\ge11$.\medskip

Thus, $e_H(T,U)\in\{n+1,n\}$. Then $e(H[T])\ge4$. Without loss of generality, $H[T]$ contains a copy of $C_r$ with vertices $w_1,\ldots,w_r$ in order. We claim that $e(H[T])=4$ and $d_H(u_1)\le2$. Suppose not. Then either $e(H[T])\ge5$ or $d_H(u_1)\ge3$. It follows that $r=4$ and $d_H(u_1)=\cdots=d_H(u_4)=2$ and $d_H(u_i)=1$ for each $i\ge5$ when $e(H[T])\ge5$;  $r=3$ and $w_4$ is adjacent to exactly one vertex of $C_r$, say $w_3$, $d_H(u_1)=3$, $d_H(u_2)=d_H(u_3)=d_H(u_4)=2$ and $d_H(u_i)=1$ for each $i\ge5$ when $d_H(u_1)\ge3$. Let $N_H(u_i)=\{w_i,w_{i+1}\}$ for each $i\in[3]$ and $\{w_1,w_{4}\}\subseteq N_H(u_4)$.   Note that $H-\{u_i,u_j\}$ contains two edge-disjoint $M_4$ for any $i,j\ge5$. By Claim 1, $E(G[U-\{u_1,\ldots,u_4\}])=\emptyset$. Hence, $e_G(u_i,\{u_1,\ldots,u_4\})\ge 1$ for each $i\ge 5$ because $\delta(G)\ge2$. Then there exist two adjacent vertices in $T$, say $w_1,w_2$, such that $d_{U-\{u_1,\ldots,u_4\}}(\{w_1,w_2\})\ge3$ because $|U|\ge11$. By symmetry, let $u_5w_1,u_6w_2,u_7w_2\in E(H)$ and $u_6u_1\in E(G)$. Note that $H-\{u_1,u_6\}$ contains two edge-disjoint $M_4$. But then $u_6u_1\notin E(G)$ by Claim 1, a contradiction, as claimed.
Thus, $r=3$ and $w_4$ is adjacent to exactly one vertex of $C_r$, say $w_3$, $d_H(u_1)=\cdots=d_H(u_5)=2$ and $d_H(u_i)=1$ for each $i\ge6$. Let $N_H(u_1)=\{w_1,w_3\}$, $N_H(u_i)=\{w_{i-1},w_i\}$ for each $2\le i\le 4$ and $N_H(u_5)=\{w_3,w_4\}$. Note that $H-\{u_i,u_j\}$ contains two edge-disjoint $M_4$ for any $i,j\ge6$. By Claim 1, $E(G[U-\{u_1,\ldots,u_5\}])=\emptyset$. Hence, $e_G(u_i,\{u_1,\ldots,u_5\})\ge 1$ for each $i\ge 6$ because $\delta(G)\ge2$. Then either $d_{U-\{u_1,\ldots,u_5\}}(\{w_1,w_2\})\ge3$ or $d_{U-\{u_1,\ldots,u_5\}}(\{w_3,w_4\})\ge4$ because $|U|\ge11$. Suppose  $d_{U-\{u_1,\dots,u_5\}}(\{w_1,w_2\})\ge3$. Let $u_6w_1,u_7w_2,u_8w_i$ and $u_8u_2\in E(G)$, where $i\in[2]$. Note that $H-\{u_8,u_2\}$ contains two edge-disjoint $M_4$. But then $u_8u_2\notin E(G)$ by Claim 1, a contradiction.  Suppose $d_{U-\{u_1,\ldots,u_5\}}(\{w_3,w_4\})\ge4$. Then there exist two vertices $u_i,u_j$ such that $\{u_iu_4,u_ju_5\}\subseteq E(G)$, where $i\ne j\ge6$. Note that $H-\{w_4,u_4,u_5,u_i,u_j\}$ contains two edge-disjoint $M_3$, say $M$ and $M'$. By Claim 2, $\alpha(u_iu_4)=\alpha(w_4u_5)$ and $\alpha(u_ju_5)=\alpha(w_4u_4)$. But then $G$ has a rainbow $M_5=M\cup\{u_iu_4,u_ju_5\}$, a contradiction.\qed\\

\noindent{\bf Acknowledgments.}\\

Yongxin Lan was partially supported by the National Natural Science Foundation of China (No. 12001154), Natural Science Foundation of Hebei Province (No. A2021202025) and the Special Funds for Jointly Building Universities of Tianjin (No. 280000307).

\end{document}